\documentclass[12pt]{article}
\usepackage{amssymb,amsmath,amsthm}
\usepackage{color}

\theoremstyle{definition}

\begin{document}

\begin{center}
{\bf
A non-existence  result due to small perturbations  

in an eigenvalue problem
}

Gelu Pa\c{s}a,

{\it 
Simion Stoilow Institute of Mathematics of Romanian Academy}

pasa.gelu@gmail.com

\end{center}

\noindent 
{\bf Abstract}. 
We consider a well-posed 
eigenvalue problem on 
$(a,0)$, depending 
on a continuous 
function $m$. The boundary conditions  in the 
points $a,0$ are 
depending on  the
eigenvalues.
We divide $(a,0)$ into small intervals 
 and   approximate the function $m$ by a simple (step) function $m_S$, constant on each small   interval.
The eigenfunctions 
corresponding to 
$m_S$ do not exist.

\vspace{0.25cm}


\hspace{1.7cm}proximation.


\vspace{0.5cm}

\noindent
{\bf 1. INTRODUCTION}.

\vspace{0.25cm}

 A well-posed 
eigenvalue problem 
is considered on the segment
$(a,0)$, 
depending on a continuous 
function $m$ as 
parameter. This problem 
is related with 
 some previous papers
concerning the linear stability of displacements in Hele-Shaw cells, see
\cite{CA-PA-1} -
\cite{PA1}.
The boundary conditions 
in the points $a,0$ are 
containing the eigenvalues - then we have some
compatibility conditions.
An existence result can be 
obtained by using the finite-difference method,
see \cite{CA-PA-1}.
We divide 
$(a,0)$
into small intervals 
 and   approximate the function $m$ by a simple (step) function $m_S$, constant on each small   interval. In this case,  the eigenfunctions 
corresponding to 
$m_S$ are known - a combination of exponentials. From this reason, the above mentioned
compatibility conditions
(related with the boundary conditions in the points
$a$ and $0$) can not be verified.
As a consequence,  we obtain 
the following  result. 
The functions
$m, m_S$   could be
are very ``close''
(in some suitable functional spaces), 
but the eigenvalues corresponding to
$m_S$ do not exist.
Therefore we get a very strong modification of 
the  eigenvalues, 
caused by the very small variation of the parameter
$m$.
A similar result  was given 
in \cite{PA3}:  the
``perturbed'' eigenvalues 
become infinity
with increasing wavenumbers.

\vspace{0.5cm}

\noindent
{\bf 2. THE ORIGINAL
 PROBLEM}
 
\vspace{0.25cm} 
 
On the segment 
$(a,0)$
we consider the
eigenvalue problem
(the lower indices $_x$ denote the partial derivatives)
$$
  -(m f_x)_x +  k^2 m f = 
  (1/ \sigma) k^2  f m_x, \quad k \geq 0,
$$
\begin{equation}\label{A}
 0< m_L \leq m(x) 
\leq m_R,   
\quad \forall x \in (a,0).
\end{equation}
The  
following conditions for large $|x|$ are supposed:
$$
f(x)= A e^{kx},
x\leq a;     \quad
f(x)= A e^{-kx},
x \geq 0;              $$
\begin{equation}\label{A1}
m(x)= 
m_L, x < a;  \quad
   m(x)= 
m_R, x > 0.       
\end{equation}
As  we specified in {\bf Introduction}, this problem is related with
some previous papers concerning the stability
of the flow displacements in Hele-Shaw cells.
In these   papers, 
the eigenfunctions $f$ are the amplitudes
of the perturbed velocities.
For this reason we consider that $f$ must decrease to zero
at large distances. 
Moreover, we suppose that
$f$ are continuous in 
$a,0$ but $f_x$ could 
have jumps in these
points.
On the other hand, the parameter
$m$ can be considered as 
 the 
viscosity of an intermediate fluid, between a displacing
and displaced 
fluids with 
constant viscosities
$m_L, m_R$,  see   the papers 
\cite{CA-PA-1} -
\cite{PA1}.
The function 
$m$ is  continuous on 
$[a,b]$.

We consider the following boundary conditions
in the points $a$ and
$b=0$: 
$$  m^-(a) f_x^-(a) - m^+(a) f_x^+(a)= 
[kE(a) / \sigma]f(a),  $$
$$  m^-(b) f_x^-(b) - m^+(b)  f_x^+(b)= 
[kE(b)/ \sigma]f(b),  $$ 
$$ 
E(a):= k 
[m^+(a)-m^-(a)]=
k[m^+(a)-m_L],
$$ 
\begin{equation}\label{LAPLACE002}
E(b):= k
[m^+(b)-m^-(b)]=
k[m_R-m^-(b)],   
\end{equation}
where 
$^-, ^+$ are denoting 
the left and right limit values.
$\sigma$ can be considered  
as the characteristic values of the problem; 
$k, f $ are  
the wavenumbers and the eigenfunctions.

In  the relation \eqref{LAPLACE002}  we have the same eigenvalue(s), therefore
$f$ must verify the compatibility 
relation
\begin{equation}
\label{LAPLACE002A}
V(a)=V(b),
\end{equation}
$$  
V(a):=\frac{kE(a)f(a)} 
{m^-(a) f_x^-(a) - 
m^+(a) f_x^+(a)},  
  $$

\vspace{0.25cm}  
  
  \begin{equation}\label{LAPLACE003} 
V(b):=\frac{kE(b)f(b)}  
{m^-(b) f_x^-(b) - 
m^+(b)  f_x^+(b)}.
\end{equation} 
It is important to note that 
\eqref{LAPLACE002A} is depending only on $f$
and is not depending 
on $\sigma$.
A quite similar 
problem 
has been studied in \cite{PA2}, where
$E(a),E(b)$  
 were polynomials of order 3 in $k$, depending also
 on the surface tensions
 between the displacing fluids.

 As in \cite{PA2},
we use 
the finite difference 
method  to get
a classical eigenvalue problem  for a matrix,
equivalent with 
\eqref{A} -
\eqref{LAPLACE002}.
We use the notation
$$ \lambda = 1/ \sigma, $$
thus $\lambda$ can be considered as eigenvalue(s).

The values
$f_x^-(a), f_x^+(b)$ are given by the above relations \eqref{A1},
because we know the values of $f$ outside the segment $(a,b)$:
\begin{equation}\label{A3}
 f_x^-(a)= kf(a),  \quad
 f_x^+(b)= -kf(b).
\end{equation}
From the relations
\eqref{LAPLACE002} we can obtain
the limit values of 
$f_x$ in  the interior ends
of the segment $(a,b)$:
$$ m^+(a) f_x^+(a)= 
m_L kf(a) 
-\lambda kE(a)f(a),  
                      $$
$$  m^-(b) f_x^-(b) = 
-m_R k f(b) + 
\lambda kE(b)f(b).    $$
We consider $(M+1)$ equidistant points $x_i$
s.t.
$$ a = x_M <x_{M-1}....<
x_1<x_0 = b =0,
\quad x_{i-1}-x_i=d, $$
and the following  
approximations for 
$f_x^+(a), f_x^-(b), f_{xx}(z)$ for  
$ \,\, z \in (a,b):$ 
$$
\frac{f_0-f_1}{d} =
f_x^+(a),  \quad
\frac{f_{M-1}-f_M}{d} =
f_x^-(b),
                      $$
\begin{equation}\label{A4}
f_{xx}(z)= \frac{f(z+d)-
2f(z)+f(z-d)}{d^2}.
\end{equation}
The last formula is obtained by using the symmetric approximation
of $f_z$ with the 
approximation step $d/2$.
Therefore, the boundary conditions become
$$
-\frac{m^+(a)}
{kdE(a)}f_{M-1}+ 
[\frac{m_Lk}{kE(a)}+
 \frac{m^+(a)}
 {kdE(a)}]f_{M}=
 \lambda f_M,
                       $$
\begin{equation}\label{A5}
[\frac{m_Rk}{kE(b)}+ 
\frac{m^-(b)}
{kdE(b)}]f_0 -
\frac{m^-(b)}{kdE(b)}f_1=
 \lambda f_0.          
\end{equation}
The discrete form of our problem is
\begin{equation}\label{A7}
 P {\bf f} = Q {\bf f},
 \quad
 {\bf f}=
 (f_0, f_1, f_2,..., f_M).
\end{equation}
In the case of 4 interior points we have
$$P=
\left( \begin{array}
{cccccc}
a_{00} & a_{01} & 0 & 0 & 0 & 0 \\
a_{10} & a_{11} & a_{12} & 0 & 0 & 0 \\
0 & a_{21} & a_{22} & a_{23} & 0 & 0 \\
0 & 0 & a_{32} & a_{33} & a_{34} & 0 \\
0 & 0 & 0 & a_{43} & a_{44} & a_{45} \\
0 & 0 & 0 & 0 & a_{54} & a_{55} \\
\end{array} \right),
$$
where
$$ a_{00}= 
[\frac{m_Rk}{kE(b)}+ 
\frac{m^-(b)}
{kdE(b)}], \quad
a_{01} = - 
\frac{m^-(b)}{kdE(b)}
                       $$
$$ a_{54}= 
-\frac{m^+(a)}{kdE(a)}
, \quad 
a_{55} = 
[\frac{m_Lk}{kE(a)}+
 \frac{m^+(a)}
 {kdE(a)}]
                       $$

    The matrix 
$P$ is  tridiagonal.
The first and the last
lines contain 
the boundary conditions
\eqref{A5}. The other
lines contain only 3
elements different from zero related to the
approximation of 
 $(mf_x)_x$ 
in $x_1,..., x_{M-1}$.

The matrix $Q$ is diagonal:
$N_{00}=N_{MM}=1$ and 
the other elements of the diagonal contain the values of $m_x$ in the interior points 
$x_1,..., x_{M-1}$. 

If 
$m_x \neq 0$ in all
interior points, then 
$Q^{-1}$   exist and 
\eqref{A7} is equivalent with
\begin{equation}\label{A8}
 Q^{-1}P {\bf f} = \lambda 
 {\bf f}.
\end{equation}
This is an eigenvalue problem for the matrix
$(Q^{-1}P )$ and  the eigenvalues 
$\lambda$ are obtained by
using the classical method.

\vspace{0.5cm}

\noindent
{\bf Remark 1.} 
If the parameter
$m$ is constant, we prove
that the eigenvalues of the problem 
 \eqref{A} -
\eqref{LAPLACE003} 
do not exist.  This fact is due to the particular form of the
boundary conditions.
When $m$ is constant,
the right hand side of
\eqref{A} is zero, thus the eigenfunctions
$f$ are known:
\begin{equation}\label{A8}
 f(x) = J e^{kx} + 
 L e^{-kx}, \quad J, L \quad
 \mbox{constant}.
\end{equation}
We consider large values of the wavenumbers $k$.
In order to avoid large values of $f$, it follows 
$$ f(x) =  J e^{kx} .   $$
We need continuity of $f$
in $a$, 
then $J=A$ where $A$ appears in the formula 
\eqref{A1}. We also need 
the continuity of $f$ in $b=0$.
Thus  we are led to   the eigenfunctions
\begin{equation}
 \label{A8B}
 f(x)= A e^{kx}, x \leq 0;
 \quad f(x)= Ae^{-kx},
 x \geq 0.
 \end{equation}

 Just now   we prove that the eigenvalue(s)
corresponding to the eigenfunctions \eqref{A8B}
do not exist. To this end,
we use  the equivalent form of the relations 
\eqref{A5}:
\begin{equation}\label{A9}
\lambda =
\frac{m^+(a)}
{kdE(a)}(-
\frac{f_{M-1}}{f_M})+ 
[\frac{m_Lk}{kE(a)}+
 \frac{m^+(a)}
 {kdE(a)}],
 \end{equation}
\begin{equation}
\label{A10}
\lambda = 
\frac{m^-(b)}{kdE(b)}
(- \frac{f_1}{f_0}) +
[\frac{m_Rk}{kE(b)}+ 
\frac{m^-(b)}
{kdE(b)}].
\end{equation}
In both above relations we have the same eigenvalue(s) $\lambda$. 
However we get
\begin{equation}
\label{A11}
(-\frac{f_{M-1}}{f_M})=- \frac{e^{k(a+d)}}
{e^{ka}}=
- e^{kd}, \quad \quad 
(- \frac{f_1}{f_0})=-
e^{-kd}.
\end{equation}
For large $k$, 
from \eqref{A9}  we get
$\lambda \rightarrow 
- \infty$
and \eqref{A10} gives us 
$\lambda \rightarrow 0 $.
As we obtain different limits for $\lambda$, 
the eigenvalue(s) do not exist.

\hfill $\square$.

\vspace{0.5cm}

\noindent
{\bf 3. THE PERTURBED
PROBLEM}

\vspace{0.25cm}

Consider the step function                            
\begin{equation}\label{VISCO-3}
 m_S(x)= m_i, \,\,
 x \in (x_{i-1}, x_i), 
\end{equation}
where
$m_i$ is constant in  each small interval. Thus we have 
\begin{equation}\label{N-0} 
-(f_S)_{xx}+k^2f_S=0, \quad
x\neq x_i, \quad 
x \in (a,0).
\end{equation}
Both   $m_S, f_S$
verify the relations 
\eqref{A1}.
The boundary conditions 
\eqref{LAPLACE002}
exist only in 
$x=a, x=b=0$ :
$$
m^-_S(a)(f_S)_x^-(a) - 
m^+_S(a)(f_S)_x^+(a) =  \frac{kE_S(a)f_S(a)}{\sigma_S},             $$
\begin{equation}\label{N-01} 
m^-_S(b)(f_S)_x^-(b) - 
m^+_S(b)(f_S)_x^+(b) =  \frac{kE_S(b)f_S(b)}{\sigma_S},             
\end{equation}
where  $E_S(a), E_S(b)$  are given by \eqref{LAPLACE002}
with $m_S$ instead of 
$m$   and $\sigma_S$ are 
the perturbed characteristic values.
As in  
{\bf Remark 1},
\eqref{N-01} 
contain  the same 
$\sigma_S$,
thus we need 
\begin{equation}\label{D1}
V_S(a)=V_S(b);
\end{equation}
$$
 V_S(a)=
 \frac{kE_S(a)f_S(a)}
{m^-_S(a)(f_S)_x^-(a) - 
m^+_S(a)(f_S)_x^+(a) }, 
                       $$

\begin{equation}\label{D2}
 V_S(b)=
\frac{kE_S(b)f_S(b)}
{m^-(b)_S(f_S)_x^-(b) - 
m^+(b)_S(f_S)_x^+(b) }.
\end{equation}
 
\vspace{0.5cm}
 
\noindent 
{\bf Remark 2.}
The problem \eqref{VISCO-3} -
\eqref{N-01}
has no solution
As in {\bf Remark 1}, 
from \eqref{N-0}
we get the following
perturbed   eigenfunctions
$$
f_S(x)= A_i e^{kx} +
B_i e^{-kx}, \quad 
x \in (x_{i-1}, x_i).
$$
Thus for large $k$ it
follows
\begin{equation}\label{D3}
 f_S(x)= A_i e^{kx}, \quad 
x \in (x_{i-1}, x_i).
\end{equation}
As 
$f$  is continuous continuous, there exists
a constant $A$ s.t.
$$  
A= A_1=A_2=...=A_N.
                        $$
Then in fact we have the same solution     \eqref{A8B}, which is continuous also in $a,0$:                     
$$ (f_S)^-(a)=f_S(a)= 
Ae^{ka};
   \quad 
   (f_S)^+(b)=f_S(b)= 
Ae^{kb}= A.            $$

We can use the same
properties of $\lambda$,
obtained in {\bf Remark 1}
and get the nonexistence of the eigenvalues, because the boundary
conditions in $a,0$ are the same.

However, we give here a different proof, based on the relation \eqref{D2}.
A jump of $f_x$ exists only  in $x_0=b=0$:
$$ (f_S)_x^-(a)=
   (f_S)_x^+(a)=
k Ae^{ka};             $$
$$
(f_S)_x^-(b)=kAe^{kb}, 
\quad 
(f_S)_x^+(b)=- kAe^{kb}
                       $$ We insert the above 
derivatives in  the relation \eqref{D2}
and it follows
$$V_S(a)= 
\frac{k(m_1-m_L)}
{m_L-m_1};   \quad       
V_S(b)= 
\frac{k(m_R-m_N)}
{m_N+ m_R}.
$$
As a consequence, we get
\begin{equation}\label{N-0A2}
V_S(a)=V_S(b) \Rightarrow 
m_R =0,
\end{equation}
which is not possible,
 because $m_R$ 
must be positive.
Thus
$V_S(a)\neq V_S(b)$.

\hfill  $\square$

\noindent
{\bf Remark 3}. 
As the relation \eqref{D1} is not  fulfilled,     
the problem 
\eqref{VISCO-3} -
\eqref{D2}   has no solution. Thus  $\sigma_S$
do not exist.
Moreover, 
the  problem \eqref{A} -
\eqref{LAPLACE003} 
can not be approximated by the  problem  
\eqref{VISCO-3} -
\eqref{D2}.

\hfill $\square$

\noindent
{\bf 4. CONCLUSION } 

\vspace{0.25cm}

The 
functions $m, m_S$ could be very ``close'' 
(in a suitable functional space) if the  number
of the jump points $x_i$
(for $m_S$) is very large. 
A very small
perturbation of a parameter can
give a strong variation of the perturbed eigenvalues - see \cite{WI}.
In our  case,
the perturbed eigenvalues do not exist, even if the
initial problem is well-posed and the ``perturbations'' are very small.

\end{document}